\documentclass[12pt, a4paper]{amsart}
\usepackage{amssymb, amsmath}
\usepackage{hyperref, upref}

\newcommand\field[1]{\mathbb{#1}}
\newcommand{\CC}{\field{C}}
\newcommand{\NN}{\field{N}}

\newcommand{\TT}{\field{T}}
\newcommand{\ZZ}{\field{Z}}
\newcommand{\Bb}{\mathcal B}

\newcommand{\Hh}{\mathcal H}

\newcommand{\Mm}{\mathcal M}

\newcommand{\Tt}{\mathcal T}

\newcommand\lsp{\mathop{\operatorname{span}}\nolimits}
\newcommand\clsp{\overline{\lsp}}

\newcommand{\spc}{\operatorname{sp}}
\newcommand{\range}{\operatorname{range}}

\newcommand{\pr}[1]{p^m_{#1}}
\newcommand{\sr}[1]{s^m_{#1}}
\newcommand{\srsp}[2]{(\sr{#1})^{#2}}
\newcommand{\srstar}[1]{\srsp{#1}{*}}
\newcommand\Cr[1]{\ensuremath{C^*_{\operatorname{min}}(#1)}}
\newcommand{\pT}[1]{p^\Tt_{#1}}
\newcommand{\sT}[1]{s^\Tt_{#1}}
\newcommand{\sTsp}[2]{(\sT{#1})^{#2}}
\newcommand{\sTstar}[1]{\sTsp{#1}{*}}
\newcommand{\Pap}[1]{P^{\Omega}_{#1}}
\newcommand{\Sap}[1]{S^{\Omega}_{#1}}
\newcommand{\Sapsp}[2]{(\Sap{#1})^{#2}}
\newcommand{\Sapstar}[1]{\Sapsp{#1}{*}}
\newcommand{\TE}{\widetilde{E}}
\newcommand{\PI}{\operatorname{PI}}
\newcommand{\Ered}{W}
\newcommand{\Erinfty}{\Omega}

\theoremstyle{plain}
\newtheorem{theorem}{Theorem}[section]
\newtheorem{cor}[theorem]{Corollary}
\newtheorem{lemma}[theorem]{Lemma}
\newtheorem{prop}[theorem]{Proposition}
\theoremstyle{definition}

\newtheorem{notation}[theorem]{Notation}
\theoremstyle{remark}
\newtheorem{rmk}[theorem]{Remark}
\newtheorem{rmks}[theorem]{Remarks}

\numberwithin{equation}{section}

\title{The co-universal \texorpdfstring{$C^*$}{C*}-algebra of a row-finite graph}
\author{Aidan Sims}
\email{asims@uow.edu.au}
\address{School of Mathematics and Applied Statistics\\
Austin Keane Building (15)\\
University of Wollongong\\
NSW 2522\\
AUSTRALIA}
\thanks{This research was supported by the Australian Research Council.}
\date{\today}

\begin{document}

\begin{abstract}
Let $E$ be a row-finite directed graph. We prove that there
exists a $C^*$-algebra $\Cr{E}$ with the following co-universal
property: given any $C^*$-algebra $B$ generated by a
Toeplitz-Cuntz-Krieger $E$-family in which all the vertex
projections are nonzero, there is a canonical homomorphism from
$B$ onto $\Cr{E}$. We also identify when a homomorphism from
$B$ to $\Cr{E}$ obtained from the co-universal property is
injective. When every loop in $E$ has an entrance, $\Cr{E}$
coincides with the graph $C^*$-algebra $C^*(E)$, but in
general, $\Cr{E}$ is a quotient of $C^*(E)$. We investigate the
properties of $\Cr{E}$ with emphasis on the utility of
co-universality as the defining property of the algebra.
\end{abstract}

\maketitle

\section{Introduction}

The aim of this paper is to initiate a study of $C^*$-algebras
defined by what we refer to as co-universal properties, and to
demonstrate the utility of such a property in investigating the
structure of the resulting $C^*$-algebra. We do this by
considering the specific example of co-universal $C^*$-algebras
associated to row-finite directed graphs.

A directed graph $E$ consists of a countable set $E^0$ of
vertices, and a countable set $E^1$ of directed edges. The
edge-directions are encoded by maps $r, s : E^1 \to E^0$: an
edge $e$ points from the vertex $s(e)$ to the vertex $r(e)$. In
this paper, we follow the edge-direction conventions of
\cite{Raeburn2005}; that is, a path in $E$ is a finite sequence
$e_1e_2...e_n$ of edges such that $s(e_i) = r(e_{i+1})$ for $1
\le i < n$.

Let $E$ be a directed graph. A \emph{Toeplitz-Cuntz-Krieger
$E$-family} in a $C^*$-algebra $B$ consists of sets $\{p_v : v
\in E^0\}$ and $\{s_e : s \in E^1\}$ of elements of $B$ such
that
\begin{enumerate}\renewcommand{\theenumi}{T\arabic{enumi}}
\item\label{it:TR1} the $p_v$ are mutually orthogonal
    projections;
\item\label{it:TR2} $s^*_e s_e = p_{s(e)}$ for all $e \in
    E^1$; and
\item\label{it:TR3} $p_v \ge \sum_{e \in F} s_e s^*_e$ for
    all $v \in E^0$ and all finite $F \subset r^{-1}(v)$.
\end{enumerate}
A Toeplitz-Cuntz-Krieger $E$-family $\{p_v : v \in E^0\}$,
$\{s_e : s \in E^1\}$ is called a \emph{Cuntz-Krieger
$E$-family} if it satisfies
\begin{enumerate}\renewcommand{\theenumi}{CK}
\item\label{it:CK} $p_v = \sum_{r(e) = v} s_e s^*_e$
    whenever $0 < |r^{-1}(v)| < \infty$.
\end{enumerate}
The graph $C^*$-algebra $C^*(E)$ is the universal $C^*$-algebra
generated by a Cuntz-Krieger $E$-family.

To see where (\ref{it:TR1})--(\ref{it:TR3}) come from, let
$E^*$ denote the path category of $E$. That is, $E^*$ consists
of all directed paths $\alpha = \alpha_1 \alpha_2 \dots
\alpha_m$ endowed with the partially defined associative
multiplication given by concatenation. There is a natural
notion of a ``left-regular" representation $\lambda$ of $E^*$
on $\ell^2(E^*)$: for a path $\alpha \in E^*$,
$\lambda(\alpha)$ is the operator on $\ell^2(E^*)$ such that
\begin{equation}\label{eq:lrr}
\lambda(\alpha) \xi_\beta =
\begin{cases}
\xi_{\alpha\beta} &\text{ if $s(\alpha) = r(\beta)$} \\
0 &\text{ otherwise.}
\end{cases}
\end{equation}
It is not hard to verify that the elements $P_v := \lambda(v)$
and $S_e := \lambda(e)$ satisfy (\ref{it:TR1})--(\ref{it:TR3}).
Indeed, it turns out that the $C^*$-algebra generated by these
$P_v$ and $S_e$ is universal for Toeplitz-Cuntz-Krieger
$E$-families.

The final relation~(\ref{it:CK}) arises if we replace the space
$E^*$ of paths in $E$ with its boundary $E^{\le\infty}$ (this
boundary consists of all the infinite paths in $E$ together
with those finite paths that originate at a vertex which
receives no edges). A formula more or less identical
to~\eqref{eq:lrr} defines a Cuntz-Krieger $E$-family
$\{P^\infty_v : v \in E^0\}$, $\{S^\infty_e : e \in E^1\}$ in
$\Bb(\ell^2(E^{\le\infty}))$. The Cuntz-Krieger uniqueness
theorem \cite[Theorem~3.1]{BPRS2000} implies that when every
loop in $E$ has an entrance, the $C^*$-algebra generated by
this Cuntz-Krieger family is universal for Cuntz-Krieger
$E$-families. When $E$ contains loops without entrances
however, universality fails. For example, if $E$ has just one
vertex and one edge, then a Cuntz-Krieger $E$-family consists
of a pair $P,S$ where $P$ is a projection and $S$ satisfies
$S^*S = P = SS^*$. Thus the universal $C^*$-algebra $C^*(E)$ is
isomorphic to $C^*(\ZZ) = C(\TT)$. However, $E^{\le\infty}$
consists of a single point, so $C^*(\{P^\infty_v, S^\infty_e\})
\cong \CC$.

The definition of $C^*(E)$ is justified, when $E$ contains
loops with no entrance, by the gauge-invariant uniqueness
theorem (originally due to an Huef and Raeburn; see
\cite[Theorem~2.3]{HR1997}), which says that $C^*(E)$ is the
unique $C^*$-algebra generated by a Cuntz-Krieger $E$-family in
which each $p_v$ is nonzero and such that there is a
\emph{gauge action} $\gamma$ of $\TT$ on $C^*(E)$ satisfying
$\gamma_z(p_v) = p_v$ and $\gamma_z(s_e) = zs_e$ for all $v \in
E^0$, $e \in E^1$ and $z \in \TT$.

Recently, Katsura developed a very natural description of this
gauge-invariant uniqueness property in terms of what we call
here a \emph{co-universal property}. In the context of graph
$C^*$-algebras, Proposition~7.14 of \cite{Katsura2007} says
that $C^*(E)$ is co-universal for gauge-equivariant
Toeplitz-Cuntz-Krieger $E$-families in which each vertex
projection is nonzero. That is, $C^*(E)$ is the unique
$C^*$-algebra such that
\begin{itemize}
\item $C^*(E)$ is generated by a Toeplitz-Cuntz-Krieger
    $E$-family $\{p_v, s_e\}$ such that each $p_v$ is
    nonzero, and $C^*(E)$ carries a gauge action; and
\item For every Toeplitz-Cuntz-Krieger $E$-family $\{q_v,
    t_e\}$ such that each $q_v$ is nonzero and such that
    there is a strongly continuous action $\beta$ of $\TT$
    on $C^*(\{q_v, t_e\})$ satisfying $\beta_z(q_v) = q_v$
    and $\beta_z(t_e) = zt_e$ for all $v \in E^0$ and $e
    \in E^1$, there is a homomorphism $\psi_{q,t} :
    C^*(\{q_v, t_e\}) \to C^*(E)$ satisfying
    $\psi_{q,t}(q_v) = p_v$ and $\psi_{q,t}(t_e) = s_e$ for
    all $v \in E^0$ and $e \in E^1$.
\end{itemize}

The question which we address in this paper is whether there
exists a co-universal $C^*$-algebra for (not necessarily
gauge-equivariant) Toeplitz-Cuntz-Krieger $E$-families in which
each vertex projection is nonzero. Our first main theorem,
Theorem~\ref{thm:Cr(E)-Existence} shows that there does indeed
exist such a $C^*$-algebra $\Cr{E}$, and identifies exactly
when a homomorphism $B \to \Cr{E}$ obtained from the
co-universal property of the latter is injective. The bulk of
Section~\ref{sec:existence} is devoted to proving this theorem.
Our key tool is Hong and Szyma\'nski's powerful description of
the primitive ideal space of the $C^*$-algebra of a directed
graph. We realise $\Tt C^*(E)$ as the universal $C^*$-algebra
of a modified graph $\TE$ to apply Hong and Szyma\'nski's
results to the Toeplitz algebra.

Our second main theorem, Theorem~\ref{thm:Cr(E)-Properties} is
a uniqueness theorem for the co-universal $C^*$-algebra. We
then proceed in the remainder of Section~\ref{sec:properties}
to demonstrate the power and utility of the defining
co-universal property of $\Cr{E}$ and of our uniqueness theorem
by obtaining the following as fairly straightforward
corollaries:
\begin{itemize}
\item a characterisation of simplicity of $\Cr{E}$;
\item a characterisation of injectivity of representations
    of $\Cr{E}$;
\item a description of $\Cr{E}$ in terms of a universal
    property, and a uniqueness theorem of Cuntz-Krieger
    type;
\item a realisation of $\Cr{E}$ as the Cuntz-Krieger
    algebra $C^*(F)$ of a modified graph $F$;
\item an isomorphism of $\Cr{E}$ with the $C^*$-subalgebra
    of $\Bb(\ell^2(E^{\le \infty}))$ generated by the
    Cuntz-Krieger $E$-family $\{P^{\le\infty}_v,
    S^{\le\infty}_e\}$ described earlier; and
\item a faithful representation of $\Cr{E}$ on a Hilbert
    space $\Hh$ such that the canonical faithful
    conditional expectation of $\Bb(\Hh)$ onto its diagonal
    subalgebra implements an expectation from $\Cr{E}$ onto
    the commutative $C^*$-subalgebra generated by the range
    projections $\{\sr{\alpha} \srstar{\alpha} : \alpha \in
    E^*\}$.
\end{itemize}

Our results deal only with row-finite graphs to simplify the
exposition. However, it seems likely that a similar analysis
applies to arbitrary graphs. Certainly Hong and Szyma\'nski's
characterisation of the primitive ideal space of a graph
$C^*$-algebra is available for arbitrary graphs. In principle
one can argue along exactly the same lines as we do in
Section~3 to obtain a co-universal $C^*$-algebra for an
arbitrary directed graph. Alternatively, the results of this
paper could be bootstrapped to the non-row-finite situation
using Drinen and Tomforde's desingularisation process
\cite{DT2005}.

\subsection*{Acknowledgements.}
The author thanks Iain Raeburn for lending a generous ear, and
for helpful conversations. The author also thanks Toke Carlsen
for illuminating discussions, and for his suggestions after a
careful reading of a preliminary draft.

\section{Preliminaries}

We use the conventions and notation for directed graphs
established in \cite{Raeburn2005}; in particular our
edge-direction convention is consistent with \cite{Raeburn2005}
rather than with, for example, \cite{BHRS2002, BPRS2000,
HS2004}.

A \emph{path} in a directed graph $E$ is a concatenation
$\lambda = \lambda_1\lambda_2 \ldots \lambda_n$ of edges
$\lambda_i \in E^1$ such that $s(\lambda_i) = r(\lambda_{i+1})$
for $i < n$; we write $r(\lambda)$ for $r(\lambda_1)$ and
$s(\lambda)$ for $s(\lambda_n)$. We denote by $E^*$ the
collection of all paths in $E$. For $v \in E^0$ we write $vE^1$
for $\{e \in E^1 : r(e) = v\}$; similarly $E^1v = \{e \in E^1 :
s(e) = v\}$.

A \emph{cycle} in $E$ is a path $\mu = \mu_1 \dots \mu_{|\mu|}$
such that $r(\mu) = s(\mu)$ and such that $s(\mu_i) \not=
s(\mu_j)$ for $1 \le i < j \le |\mu|$. Given a cycle $\mu$ in
$E$, we write $[\mu]$ for the set
\[
[\mu] = \{\mu,\; \mu_2\mu_3\cdots\mu_{|\mu|}\mu_1,\; \dots,\; \mu_{|\mu|}
\mu_1\cdots\mu_{n-1}\}
\]
of cyclic permutations of $\mu$. We write $[\mu]^0$ for the set
$\{s(\mu_i) : 1 \le i \le |\mu|\} \subset E^0$, and $[\mu]^1$
for the set $\{\mu_i : 1 \le i \le |\mu|\} \subset E^1$. Given
a cycle $\mu$ in $E$ and a subset $M$ of $E^0$ containing
$[\mu]^0$, we say that $\mu$ has \emph{no entrance in M} if
$r(e) = r(\mu_i)$ and $s(e) \in M$ implies $e = \mu_i$ for all
$1 \le i \le |\mu|$. We denote by $C(E)$ the set $\{[\mu] :
\mu$ is a cycle with no entrance in $E^0\}$. By $C(E)^1$ we
mean $\bigcup_{C \in C(E)} C^1$, and by $C(E)^0$ we mean
$\bigcup_{C \in C(E)} C^0$.

A \emph{cutting set} for a directed graph $E$ is a subset $X$
of $C(E)^1$ such that for each $C \in C(E)$, $X \cap C^1$ is a
singleton. Given a cutting set $X$ for $E$, for each $x \in X$,
we write $\mu(x)$ for the unique cycle in $E$ such that $r(\mu)
= r(x)$, and let $\lambda(x) = \mu(x)_2\mu(x)_3 \dots
\mu(x)_{|\mu(x)|}$; so $\mu(x) = x\lambda(x)$ for all $x \in
X$, and $C(E) = \{[\mu(x)] : x \in X\}$.

\section{Existence of the co-universal \texorpdfstring{$C^*$}{C*}-algebra}\label{sec:existence}

Our main theorem asserts that every row-finite directed graph
admits a co-universal $C^*$-algebra and identifies when a
homomorphism obtained from the co-universal property is
injective.

\begin{theorem}\label{thm:Cr(E)-Existence}
Let $E$ be a row-finite directed graph.
\begin{enumerate}
\item\label{it:existence} There exists a $C^*$-algebra
    $\Cr{E}$ which is co-universal for
    Toeplitz-Cuntz-Krieger $E$-families of nonzero partial
    isometries in the sense that $\Cr{E}$ is generated by a
    Toeplitz-Cuntz-Krieger $E$-family $\{P_v : v \in E^0\},
    \{S_e : e \in E^1\}$ with the following two properties.
    \begin{enumerate}
    \item\label{it:generators} The vertex projections
        $\{P_v : v \in E^0\}$ are all nonzero.
    \item\label{it:co-universal} Given any
        Toeplitz-Cuntz-Krieger $E$-family $\{q_v : v
        \in E^0\}, \{t_e : e \in E^1\}$ such that each
        $q_v \not= 0$ and given any cutting set $X$ for
        $E$, there is a function $\kappa : X \to \TT$
        and a homomorphism $\psi_{q,t} : C^*(\{q_v, t_e
        : v \in E^0, e \in E^1\}) \to \Cr{E}$
        satisfying $\psi_{q,t}(q_v) = P_v$ for all $v
        \in E^0$, $\psi_{q,t}(t_e) = S_e$ for all $e
        \in E^1 \setminus X$, and $\psi_{q,t}(t_x) =
        \kappa(x)S_x$ for all $x \in X$.
    \end{enumerate}
\item\label{it:injectivity} Given a Toeplitz-Cuntz-Krieger
    $E$-family $\{q_v : v \in E^0\}$, $\{t_e : e \in E^1\}$
    with each $q_v$ nonzero, the homomorphism $\psi_{q,t} :
    B \to \Cr{E}$ obtained from~(\ref{it:co-universal}) is
    an isomorphism if and only if for each cycle $\mu$ with
    no entrance in $E$, the partial isometry $t_\mu$ is a
    scalar multiple of $q_{r(\mu)}$.
\end{enumerate}
\end{theorem}

\begin{rmks}
It is convenient in practise to work with cutting sets $X$ and
functions from $X$ to $\TT$ as in
Theorem~\ref{thm:Cr(E)-Existence}(\ref{it:co-universal}).
However, property~(\ref{it:co-universal}) can also be
reformulated without respect to cutting sets. Indeed:
\begin{enumerate}
\item The asymmetry arising from the choice of a cutting
    set $X$ in
    Theorem~\ref{thm:Cr(E)-Existence}(\ref{it:co-universal})
    can be avoided. The property could be reformulated
    equivalently as follows: \emph{given a
    Toeplitz-Cuntz-Krieger $E$-family $\{q_v : v \in E^0\},
    \{t_e : e \in E^1\}$ such that each $q_v \not= 0$,
    there is a function $\rho : C(E)^1 \to \TT$ and a
    homomorphism $\psi_{q,t} : C^*(\{q_v, t_e : v \in E^0,
    e \in E^1\}) \to \Cr{E}$ satisfying $\psi_{q,t}(q_v) =
    P_v$ for all $v \in E^0$, $\psi_{q,t}(t_e) = S_e$ for
    all $e \in E^1 \setminus C(E)^1$, and $\psi_{q,t}(t_e)
    = \rho(e) S_e$ for all $e \in C(E)^1$.} One can prove
    that an algebra satisfying this modified
    condition~(\ref{it:co-universal}) exists  using exactly
    the same argument as for the current theorem after
    making the appropriate modification to
    Lemma~\ref{lem:Ikappas same}. That the resulting
    algebra coincides with $\Cr{E}$ follows from
    applications of the co-universal properties of the two
    algebras.
\item Fix a row-finite graph $E$ with no sources and a
    function $\kappa : C(E) \to \TT$. Let $\{q_v : v \in
    E^0\}, \{t_e : e \in E^1\}$ be a Toeplitz-Cuntz-Krieger
    $E$-family such that each $q_v \not= 0$. Then there is
    a homomorphism as in
    Theorem~\ref{thm:Cr(E)-Existence}(\ref{it:co-universal})
    with respect to the fixed function $\kappa$ for some
    cutting set $X$ if and only if there is such a
    homomorphism for every cutting set $X$. One can see
    this by following the argument of Lemma~\ref{lem:what's
    the big ideal} below to see that $\kappa$ does not
    depend on $X$.
\item Given a Toeplitz-Cuntz-Krieger $E$-family $\{q_v : v
    \in E^0\}, \{t_e : e \in E^1\}$ such that each $q_v
    \not= 0$ and a cutting set $X$, the functions $\kappa :
    X \to \TT$ which can arise in
    Theorem~\ref{thm:Cr(E)-Existence}(\ref{it:co-universal})
    are precisely those for which $\kappa([\mu])$ belongs
    to the spectrum $\spc_{q_v C^*(\{q_v, t_e : v \in E^0,
    e \in E^1\}) q_v}(t_\mu)$ for each cycle $\mu$ without
    an entrance in $E$. To see this, one uses Hong and
    Szyma\'nski's theorems to show that in the first
    paragraph of the proof of Lemma~\ref{lem:what's the big
    ideal}, the complex numbers $z$ which can arise are
    precisely the elements of the spectrum of the unitary
    $s_{\alpha(\mu)} + I$ in the corner
    $(p_{\alpha(r(\mu))} + I) \big(C^*(\TE)/I\big)
    (p_{\alpha(r(\mu))} + I)$.
\end{enumerate}
\end{rmks}

\begin{cor}
Let $E$ be a row-finite directed graph in which every cycle has
an entrance.  Then $C^*(E) \cong \Cr{E}$. In particular, if
$\{q_v : v \in E^0\}$, $\{t_e : e \in E^1\}$ is a
Toeplitz-Cuntz-Krieger $E$-family in a $C^*$-algebra $B$ such
that each $q_v$ is nonzero, then there is a homomorphism
$\psi_{q,t} : C^*(\{q_v, t_e : v\in E^0, e \in E^1\}) \to
C^*(E)$ such that $\psi_{q,t}(q_v) = p_v$ for all $v \in E^0$
and $\psi_{q,t}(t_e) = s_e$ for all $e \in E^1$.
\end{cor}
\begin{proof}
For the first statement, observe that the co-universal property
of $\Cr{E}$ induces a surjective homomorphism $\psi_{p,s} :
C^*(E) \to \Cr{E}$. Since every cycle in $E$ has an entrance,
the condition in
Theorem~\ref{thm:Cr(E)-Existence}(\ref{it:injectivity}) is
trivially satisfied, and it follows that $\psi_{p,s}$ is an
isomorphism.

Since every cycle in $E$ has an entrance, a cutting set for $E$
has no elements. Hence the second statement is just a
re-statement of
Theorem~\ref{thm:Cr(E)-Existence}(\ref{it:co-universal}).
\end{proof}

The remainder of this section will be devoted to proving
Theorem~\ref{thm:Cr(E)-Existence}. Our key technical tool in
proving Theorem~\ref{thm:Cr(E)-Existence} will be Hong and
Szyma\'nski's description of the primitive ideal space of a
graph $C^*$-algebra. To do this, we first realise the Toeplitz
algebra of $C^*(E)$ as a graph algebra in its own right. This
construction is known, but we have found it difficult to pin
down in the literature.

\begin{notation}\label{ntn:TE}
Let $E$ be a directed graph. Define a directed graph $\TE$ as
follows:
\begin{gather*}
\TE^0 = \{\alpha(v) : v \in E^0\} \sqcup \{\beta(v) : v \in
E^0, 0 < |vE^1| < \infty\} \\
\TE^1 = \{\alpha(e) : e \in E^1\} \sqcup \{\beta(e) : e \in
E^1, 0 < |s(e) E^1| < \infty\} \\
r(\alpha(e)) = r(\beta(e)) = \alpha(r(e)),\\
s(\alpha(e)) = \alpha(s(e)),\text{ and } s(\beta(e)) =
\beta(s(e)).
\end{gather*}
For $\lambda \in E^*$ with $|\lambda| \ge 2$, we define
$\alpha(\lambda) :=
\alpha(\lambda_1)\dots\alpha(\lambda_{|\lambda|})$. Since $E$
is row-finite, $\TE$ is also row-finite.
\end{notation}

\begin{lemma}\label{lem:TE alg}
Let $E$ be a directed graph and let $\TE$ be as in
Notation~\ref{ntn:TE}. For $v \in E^0$ and $e \in E^1$, let
\begin{align*}
q_v &:=
    \begin{cases}
        p_{\alpha(v)} + p_{\beta(v)} & \text{ if $0 < |vE^1| < \infty$} \\
        p_{\alpha(v)} & \text{ otherwise,}
    \end{cases} \\
\intertext{and}
t_e &:=
    \begin{cases}
        s_{\alpha(e)} + s_{\beta(e)} & \text{ if $0 < |s(e)E^1| < \infty$}\\
        s_{\alpha(e)} & \text{ otherwise.}
    \end{cases}
\end{align*}
Then there is an isomorphism $\phi : \Tt C^*(E) \to C^*(\TE)$
satisfying $\phi(\pT{v} = q_v$ and $\phi(\sT{e}) = t_e$ for all
$v \in E^0$ and $e \in E^1$.

\end{lemma}
\begin{proof}
Routine calculations show that $\{q_v : v \in E^0\}$, $\{t_e :
e \in E^1\}$ is a Toeplitz-Cuntz-Krieger $E$-family in
$C^*(\TE)$. The universal property of $\Tt C^*(E)$ therefore
implies that there is a homomorphism $\phi : \Tt C^*(E) \to
C^*(\TE)$ satisfying $\phi(\pT{v}) = q_v$ for all $v \in E^0$
and $\phi(\sT{e}) = t_e$ for all $e \in E^1$.

To see that $\phi$ is surjective, fix $v \in \TE^0$ and $E \in
\TE^1$. To see that $p_v\in \range(\phi)$, we consider three
cases: (a) $v = \alpha(w)$ for some $w$ with $wE^1$ either
empty or infinite; (b) $v = \alpha(w)$ for some $w$ with $0 <
|wE^1| < \infty$; or (c) $v = \beta(w)$ for some $w$ with $0 <
|wE^1| < \infty$. In case~(a), we have $p_v = p_{\alpha(w)} =
\phi(\pT{w})$ by definition. In case~(b), the set $v\TE^1 =
\{\alpha(e), \beta(e) : e \in wE^1\}$ is nonempty and finite.
Hence the Cuntz-Krieger relation in $C^*(\TE)$ ensures that
\begin{equation}\label{eq:CK for Toeplitz}
p_v = p_{\alpha(w)}
    = \sum_{e \in vE^1} s_{\alpha(e)} s^*_{\alpha(e)} + s_{\beta(e)} s^*_{\beta(e)}
    = \sum_{e \in vE^1} t_e t^*_e \in \range(\phi).
\end{equation}
In case~(c), we have $p_v = q_v - p_{\alpha(w)} \in
\range(\phi)$ by case~(b). Now to see that $s_e \in
\range(\phi)$, observe that if $e = \alpha(f)$ for some $f \in
E^1$, then $s_e = s_{\alpha(f)} = \phi(\sT{f}) p_{\alpha(s(f))}
\in \range(\phi)$, and if $e = \beta(f)$, then $s_e =
s_{\beta(f)} = \phi(\sT{f}) p_{\beta(s(f))} \in \range(\phi)$
also.

To finish the proof, observe that if $0 < |vE^1| < \infty$,
then $q_v - \sum_{r(e) = v} t_e t^*_e = p_{\beta(v)} \not= 0$.
Since the $t_e$ are all nonzero and have mutually orthogonal
ranges, it follows that for each $v \in E^0$ and each finite
subset $F$ of $vE^1$, we have $q_v - \sum_{e \in F} t_e t^*_e
\not= 0$. Thus the uniqueness
theorem~\cite[Theorem~4.1]{FR1999} for $\Tt C^*(E)$ implies
that $\phi$ is injective.
\end{proof}

\begin{notation}\label{ntn:Ikappa}
Let $E$ be a directed graph.
\begin{enumerate}
\item For $v \in E^0$ such that $0 < |vE^1| < \infty$, we
    define $\Delta_v := \pT{v} - \sum_{r(e) = v}
    \sT{e}\sTstar{e} \in \Tt C^*(E)$.
\item Given a function $\kappa : C(E) \to \TT$, we denote
    by $I^\kappa$ the ideal of $\Tt C^*(E)$ generated by
    $\{\Delta_v : v \in E^0\} \cup \{\kappa(C) \pT{r(\mu)}
    - \sT{\alpha(\mu)} : C \in C(E), \mu \in C\}$.
\end{enumerate}
\end{notation}

\begin{lemma}\label{lem:Ikappas same}
Let $E$ be a directed graph. Let $\kappa$ be a function from
$C(E) \to \TT$, and let $1 : C(E) \to \TT$ denote the constant
function $1(C) = 1$ for all $C \in C(E)$. Fix a cutting set $X$
for $E$, and for each $x \in X$, let $C(x)$ be the unique
element of $C(E)$ such that $x \in C(x)^1$. Then there is an
isomorphism $\widetilde{\tau_{\overline{\kappa}}} : \Tt C^*(E)
/ I^1 \to \Tt C^*(E) / I^\kappa$ satisfying
\[\begin{array}{rcll}
\widetilde{\tau_{\overline{\kappa}}}(\pT{v}+ I^1) &=& p_v + I^\kappa &\qquad\text{for all $v \in E^0$} \\
\widetilde{\tau_{\overline{\kappa}}}(\sT{e} + I^1) &=& s_e + I^\kappa &\qquad\text{for all $e \in E^1 \setminus X$, and} \\
\widetilde{\tau_{\overline{\kappa}}}(\sT{x} + I^1) &=& \overline{\kappa(C(x))}s_x + I^\kappa &\qquad\text{for all $x \in X$.}
\end{array}\]
\end{lemma}
\begin{proof}
By the universal property of $\Tt C^*(E)$, there is an action
$\tau$ of $\TT^{C(E)}$ on $\Tt C^*(E)$ such that for $\rho \in
\TT^{C(E)}$, we have
\[\begin{array}{rcll}
{\tau_\rho}(\pT{v}) &=& \pT{v} &\qquad\text{for all $v \in E^0$} \\
{\tau_\rho}(\sT{e}) &=& \sT{e} &\qquad\text{for all $e \in E^1 \setminus X$, and} \\
{\tau_\rho}(\sT{x}) &=& \rho(C(x)) \sT{x} &\qquad\text{for all $x \in X$.}
\end{array}\]
By definition of $I^1$ and $I^\kappa$ and of the action $\tau$,
we have $\tau_{\overline{\kappa}}(I^1) = I^\kappa$. Hence
$\tau_{\overline{\kappa}}$ determines an isomorphism
\[
\widetilde{\tau_{\overline{\kappa}}} : \Tt C^*(E)/I^1 \to \tau_{\overline{\kappa}}(\Tt C^*(E)) / I^\kappa = \Tt C^*(E) / I^\kappa
\]
satisfying $\widetilde{\tau_{\overline{\kappa}}}(a + I^1) =
\tau_{\overline{\kappa}}(a) + I^\kappa$ for all $a \in \Tt
C^*(E)$.
\end{proof}

\begin{lemma}\label{lem:Ikappa vert proj free}
Let $E$ be a directed graph. Fix a function $\kappa : C(E) \to
\TT$. Then $\sT{v} \not\in I^\kappa$ for all $v \in E^0$.
\end{lemma}

To prove this lemma, we need a little notation.

\begin{notation}\label{ntn:bdryCKfam}
Given a directed graph $E$, we denote by $E^{\le\infty}$ the
collection $E^\infty \cup \{\alpha \in E^* : s(\alpha)E^1 =
\emptyset\}$. There is a Cuntz-Krieger $E$-family in
$\Bb(\ell^2(E^{\le\infty}))$ determined by
\[
P^\infty_v \xi_x
    = \begin{cases}
        \xi_x &\text{ if $r(x) = v$} \\
        0 &\text{ otherwise,}
    \end{cases}
\]
and
\[
S^\infty_e \xi_x
    = \begin{cases}
        \xi_{ex} &\text{ if $r(x) = s(e)$} \\
        0 &\text{ otherwise.}
    \end{cases}
\]
\end{notation}

If $\mu$ is a cycle with no entrance in $E$, then $r(\mu)E^{\le
\infty} = \{\mu^\infty\}$, so $S^\infty_{\mu} =
P^\infty_{r(\mu)}$.

\begin{proof}[Proof of Lemma~\ref{lem:Ikappa vert proj free}]
By Lemma~\ref{lem:Ikappas same}, it suffices to show that $I^1$
contains no vertex projections. Let $\pi_{P^\infty, S^\infty} :
\Tt C^*(E) \to \Bb(\ell^2(E^{\le\infty}))$ be the
representation obtained from then universal property of $\Tt
C^*(E)$ applied to the Cuntz-Krieger family of
Notation~\ref{ntn:bdryCKfam}. Then $\ker(\pi_{P^\infty,
S^\infty})$ contains all the generators of $I^1$, so $I^1
\subset \ker(\pi_{P^\infty, S^\infty})$. Moreover,
$\ker(\pi_{P^\infty, S^\infty})$ contains no vertex projections
because each vertex of $E$ is the range of at least one $x \in
E^{\le \infty}$.
\end{proof}

From this point onward we make the simplifying assumption that
our graphs are row-finite. Though there is no obvious
obstruction to our analysis without this restriction, the added
generality would complicate the details of our arguments. In
any case, if the added generality should prove useful, it
should not be difficult to bootstrap our results to the
non-row-finite setting by means of the Drinen-Tomforde
desingularisation process applied to the graph $\TE$ of
Notation~\ref{ntn:TE}.

\begin{prop}\label{prop:Toeplitz big ideal}
Let $I$ be an ideal of $\Tt C^*(E)$ such that $\pT{v} \not\in
I$ for all $v \in E^0$. There is a function $\kappa : C(E) \to
\TT$ such that $I \subset I^\kappa$.
\end{prop}

To prove the proposition, we first establish our key technical
lemma. This lemma is implicit in Hong and Szyma\'nski's
description \cite{HS2004} of the primitive ideal space of
$C^*(\TE)$, but it takes a little work to tease a proof of the
statement out of their two main theorems.

\begin{lemma}\label{lem:what's the big ideal}
Let $\TE$ be the directed graph of Notation~\ref{ntn:TE}. Let
$I$ be an ideal of $C^*(\TE)$ such that $p_{\alpha(v)} \not\in
I$ for all $v \in E^0$. There is a function $\kappa : C(E) \to
\TT$ such that $I$ is contained in the ideal $J^\kappa$ of
$C^*(\TE)$ generated by $\{p_{\beta(v)} : v \in E^0\}$ and
$\{\kappa(C) p_{\alpha(r(\mu))} - s_{\alpha(\mu)} : C \in C(E),
\mu \in C\}$.
\end{lemma}

Before proving the lemma, we summarise some notation and
results of \cite{HS2004} as they apply to the row-finite
directed graph $\TE$ in the situation of Lemma~\ref{lem:what's
the big ideal}. A \emph{maximal tail} of $\TE$ is a subset $M$
of $\TE^0$ such that
\begin{itemize}
\item[(MT1)] $w \in M$ and $v \TE^* w \not= \emptyset$
    imply $v \in M$;
\item[(MT2)] if $v \in M$ and $v\TE^1 \not= \emptyset$,
    then there exists $e \in v\TE^1$ such that $s(e) \in
    M$; and
\item[(MT3)] if $u, v \in M$, then there exists $w \in M$
    such that $u \TE^* w \not= \emptyset$ and $v \TE^* w
    \not= \emptyset$.
\end{itemize}

We denote by $\Mm_\gamma(\TE)$ the collection of maximal tails
$M$ of $\TE$ such that every cycle $\mu$ satisfying $[\mu]^0
\subset M$ has an entrance in $M$. We denote by $\Mm_\tau(\TE)$
the collection of maximal tails $M$ of $\TE$ such that there is
a cycle $\mu$ in $\TE$ for which $[\mu]^0 \subset M$ but $\mu$
has no entrance in $M$. Since each $\beta(v)$ is a source in
$\TE$, the cycles in $\TE$ are the paths of the form
$\alpha(\mu)$ where $\mu$ is a cycle in $E$. Moreover
$\alpha(\mu) \in C(\TE)$ if and only if $\mu \in C(E)$. Thus if
$M \in M_\tau(\TE)$, then there is a unique $C_M \in C(E)$ such
that $\alpha(C_M^0) \subset M$ and $\alpha(\mu)$ has no
entrance in $M$ for each $\mu \in C_M$. We may recover $M$ from
$C_M$:
\[
M = \{\alpha(v) : v \in E^0, v E^* C_M^0 \not=
\emptyset\}.
\]

The gauge-invariant primitive ideals of $C^*(\TE)$ are indexed
by $\Mm_\gamma(\TE)$; specifically, $M \in \Mm_\gamma(\TE)$
corresponds to the primitive ideal $\PI^\gamma_M$ generated by
$\{p_w : w \in \TE^0 \setminus M\}$. The non-gauge-invariant
primitive ideals of $C^*(\TE)$ are indexed by $\Mm_\tau(\TE)
\times \TT$; specifically, the pair $(M,z)$ corresponds to the
primitive ideal $\PI^\tau_{M,z}$ generated by $\{p_w : w \in
E^0 \setminus M\} \cup \{zp_{r(\mu)} - s_\mu\}$ for any $\mu
\in C_M$ (the ideal does not depend on the choice of $\mu \in
C_M$). Corollary~3.5 of~\cite{HS2004} describes the closed
subsets of the primitive ideal space of $C^*(\TE)$ in terms of
subsets of $\Mm_\gamma(\TE) \sqcup \Mm_\tau(\TE) \times \TT$.

\begin{proof}[Proof of Lemma~\ref{lem:what's the big ideal}]
We begin by constructing the function $\kappa$. Fix $C \in
C(E)$. Since $p_{\alpha(v)} \not\in I$ for each $v \in C^0$, we
may fix a primitive ideal $J^{C}$ of $C^*(\TE)$ such that $I
\subset J^{C}$ and $p_{\alpha(v)} \not\in J^{C}$ for $v \in
C^0$. By \cite[Corollary~2.12]{HS2004}, we have either $J^{C} =
\PI^\gamma_M$ for some $M \in \Mm_\gamma(\TE)$, or $J^{C} =
\PI^\tau_{M,z}$ for some $M \in \Mm_\tau(\TE)$ and $z \in \TT$.
Since $p_{\alpha(v)} \not\in J^{C}$ for $v \in C^0$, we must
have $\alpha(C^0) \subset M$, and then the maximal tail
condition forces
\[
M = M_C := \{\alpha(v) : v \in E^0, v E^* C^0
\not= \emptyset\} \in M_\tau(\TE).
\]
Hence $J^{C} = \PI^\tau_{M,z}$ for some $z \in \TT$; we set
$\kappa(C) := z$.

For $C \in C(E)$ and $v \in C^0$, let $J^v :=
\PI^\tau_{M_C,\kappa(C)}$. Since $\beta(v) \not \in M_C$ for
all $v \in E^0$, \cite[Lemma~2.8]{HS2004} implies that
$p_{\beta(v)} \in J^v$ for all $v \in E^0$, and our definition
of $J^v$ ensures that $p_{\alpha(v)} \not \in J^v$.

We claim that the assignment $v \mapsto J^v$ of the preceding
paragraph extends to a function $v \mapsto J^v$ from $E^0$ to
the primitive ideal space of $C^*(\TE)$ such that for every $v
\in E^0$,
\begin{enumerate}
\item $I \subset J^v$;
\item\label{it:no p(alpha)s} $p_{\alpha(v)} \not\in J^v$;
\item $p_{\beta(w)} \in J^v$ for all $w \in E^0$; and
\item\label{it:kappa(C) or gi} either $J^v = \PI^\tau_{M_C,
    \kappa(C)}$ for some $C \in C(E)$, or $J^v$ is
    gauge-invariant.
\end{enumerate}
To prove the claim, fix $v \in E^0$. If $v E^* C^0
\not=\emptyset$ for some $C \in C(E)$, then $J^v :=
\PI^\tau_{M_C, \kappa(C)}$ has the desired properties. So
suppose that $v E^* C^0 = \emptyset$ for all $C \in C(E)$. Then
there exists $x^v \in vE^{\le\infty}$ such that $x^v$ does not
have the form $\lambda\mu^\infty$ for any $\lambda \in E^*$ and
cycle $\mu \in E$. The set
\[
M^v := \{\alpha(w) : w \in E^0, w E^* x^v(n) \not=\emptyset\text{ for some }n \in \NN\}
\]
is a maximal tail of $\TE$ which contains $\alpha(v)$ and does
not contain $\beta(w)$ for any $w \in E^0$. By construction of
$x^v$, every cycle in $M^v$ has an entrance in $M^v$, and so
$J^v := \PI^\gamma_{M^v}$ satisfies (\ref{it:no
p(alpha)s})--(\ref{it:kappa(C) or gi}). It therefore suffices
to show that $I \subset \PI^\gamma_{M^v}$. For each $n \in
\NN$, we have $p_{\alpha(x^v(n))} \not\in I$, so there is a
primitive ideal $J$ of $C^*(\TE)$ containing $I$ and not
containing $p_{\alpha(x^v(n))}$. The set $M_n = \{w \in \TE^0 :
p_w \not\in J\}$ is a maximal tail of $\TE$, so either $M_n \in
\Mm_\gamma(\TE)$ and $J = \PI^\gamma_{M_n}$, or $M_n \in
\Mm_\tau(\TE)$ and $J = \PI^\tau_{M_n, z}$ for some $z \in
\TT$. By definition, $x^v(n) \in M_n$, and then~(MT1) forces
$\alpha(w) \in M_n$ for all $w \in E^0$ such that $w E^* x^v(m)
\not= \emptyset$ for some $m \le n$. Hence $M^v \subset \cup_{n
\in \NN} M_n$. Hence parts (1)~and~(3) of
\cite[Corollary~3.5]{HS2004} imply that $J^v$ belongs to the
closure of the set of primitive ideals of $C^*(\TE)$ which
contain $I$, and hence contains $I$ itself. This proves the
claim.

Observe that $I \subset \bigcap_{v \in E^0} J^v$. To prove the
result, it therefore suffices to show that $\bigcap_{v \in E^0}
J^v$ is generated by $\{p_{\beta(v)} : v \in E^0\}$ and
$\{\kappa(C) p_{\alpha(r(\mu))} - s_{\alpha(\mu)} : C \in C(E),
\mu \in C\}$.

For this, first observe that for $v \in E^0$, we have
$p_{\beta(v)} \in \bigcap_{v \in E^0} J^v$ because
$p_{\beta(v)}$ belongs to each $J^v$. Fix $C \in C(E)$ and $\mu
\in C$. Since the cycle without an entrance belonging to a
given maximal tail of $E$ is unique, for each $v \in E^0$, we
have either $J^v = \PI^\tau_{M_C, \kappa(C)}$, or $p_{r(\mu)}
\in J^v$. In particular, $\kappa(C) p_{\alpha(r(\mu))} -
s_{\alpha(\mu)} \in J^v$ for each $v \in E^0$, so $\kappa(C)
p_{\alpha(r(\mu))} - s_{\alpha(\mu)} \in \bigcap_{v \in E^0}
J^v$. Hence $J^\kappa \subset \bigcap_{v \in E^0} J^v$, and it
remains to establish the reverse inclusion.

Fix a primitive ideal $J$ of $C^*(\TE)$ which contains all the
generators of $J^\kappa$. It suffices to show that $\bigcap_{v
\in E^0} J^v \subset J$. Under the bijection between primitive
ideals of $C^*(\TE)$ and elements of $\Mm_\gamma(\TE) \sqcup
\Mm_\tau(\TE) \times \TT$, the collection $\{J_v : v \in E^0\}$
is sent to $\{M^v : v E^* C^0 = \emptyset\text{ for all }C \in
C(E)\} \sqcup \{(M_C, \kappa(C)) : C \in C(E)\}$. Since each
$J^v$ trivially contains $\bigcap_{v \in E^0} J^v$, it
therefore suffices to show that the element $N_J$ of
$\Mm_\gamma(\TE) \sqcup \Mm_{\tau}(\TE) \times \TT$
corresponding to $J$ satisfies
\begin{equation}\label{eq:MJ in closure}
N_J \in \overline{\{M^v : v E^*
C^0 = \emptyset\text{ for all }C \in C(E)\}} \sqcup \overline{\{(M_C,
\kappa(C)) : C \in C(E)\}}.
\end{equation}
Let $M_J := \{v \in \TE^0 : p_v \not\in J\}$. Then either $J$
is gauge-invariant and $N_J = M_J$, or $J$ is not
gauge-invariant, and $N_J = (M_J, z)$ for some $z \in \TT$.
Observe that
\begin{equation}\label{eq:MJ in alphas}
\begin{split}
M_J &\subset \{\alpha(v) : v \in E^0\} \\
    &= \Big(\bigcup\big\{M^v : v E^* C^0 = \emptyset\text{ for all }C \in C(E)\big\}\Big) \\
    &\hskip4cm\cup \Big(\bigcup\big\{(M_C, \kappa(C)) : C \in C(E)\big\}\Big).
\end{split}
\end{equation}
We now consider three cases.

Case~1: $J$ is gauge-invariant. Then $M_J \in M_\gamma(\TE)$,
and $N_J = M_J$. In this case, \eqref{eq:MJ in alphas} together
with parts (1)~and~(3) of \cite[Corollary~3.5]{HS2004}
give~\eqref{eq:MJ in closure}.

Case~2: $J$ is not gauge-invariant, and $M_J \not= M_C$ for all
$C \in C(E)$. We have $N_J = (M_J, z)$ for some $z \in \TT$,
and since $M_J \not= M_C$ for all $C \in C(E)$, it follows that
$N_J$ does not belong to the subset $\{(M_C, \kappa(C)) : C \in
C(E)\}_{\rm min}$ of $\{(M_C, \kappa(C)) : C \in C(E)\}$
defined on page~58 of \cite{HS2004}. Hence parts (2)~and~(4ii)
of \cite[Corollary~3.5]{HS2004} imply~\eqref{eq:MJ in closure}.

Case~3: $M_J = M_C$ for some $C \in C(E)$. Fix $\mu \in C$.
Since $J$ contains $\kappa(C)p_{\alpha(r(\mu))} -
s_{\alpha(\mu)}$, we have $N_J = (M_C, \kappa(C))$ and then
part~(4iii) of \cite[Corollary~3.5]{HS2004}
implies~\eqref{eq:MJ in closure}. This completes the proof.
\end{proof}

\begin{proof}[Proof of Proposition~\ref{prop:Toeplitz big
ideal}] Let $\phi : \Tt C^*(E) \to C^*(\TE)$ be the isomorphism
of Lemma~\ref{lem:TE alg}. Observe that by~\eqref{eq:CK for
Toeplitz}, we have $\phi(\Delta_v) = p_{\beta(v)}$ for all $v
\in E^0$ such that $vE^1 \not=\emptyset$. We claim that
$p_{\alpha(v)} \not\in \phi(I)$ for all $v \in E^0$. To see
this, first suppose that $vE^1 = \emptyset$. Then
$p_{\alpha(v)} = \phi(\pT{v}) \not\in \phi(I)$ by assumption.
Now suppose that $vE^1 \not= \emptyset$, say $r(e) = v$. Then
\[
s_{\alpha(e)}^* p_{\alpha(v)} s_{\alpha(e)} + s_{\beta(e)}^* p_{\alpha(v)} s_{\beta(e)}
    = p_{\alpha(s(e))} + p_{\beta(s(e))} = \phi(\pT{s(e)}) \not\in \phi(I)
\]
by assumption. This forces $p_{\alpha(v)} \not\in \phi(I)$.

Lemma~\ref{lem:what's the big ideal} therefore applies to the
ideal $\phi(I)$ of $C^*(\TE)$. Let $\kappa : C(E) \to \TT$ and
$J^\kappa \lhd C^*(\TE)$ be the resulting function and ideal.
Then $I^\kappa := \phi^{-1}(J^\kappa)$ is generated by
$\{\Delta_v : v \in E^0\}$ and $\{\kappa(C) \pT{r(\mu)} -
\sT{\alpha(\mu)} : C \in C(E), \mu \in C\}$ by definition of
$\phi$, and contains $I$ by construction.
\end{proof}

We are now ready to prove our main theorem.

\begin{proof}[Proof of Theorem~\ref{thm:Cr(E)-Existence}]
With $I^1 \lhd \Tt C^*(E)$ defined as in Notation~3.6 and
Lemma~3.7, define $\Cr{E} := \Tt C^*(E) / I^1$. For $v \in E^0$
and $e \in E^1$, let $P_v := \pT{v} + I^1$ and $S_e := \sT{v} +
I^1$. Then $\{P_v : v \in E^0\}$, $\{S_e : e \in E^1\}$ is a
Cuntz-Krieger $E$-family which generates $\Cr{E}$. The $P_v$
are all nonzero by Lemma~\ref{lem:Ikappa vert proj free}.

Now let $\{q_v : v \in E^0\}, \{t_e : e \in E^1\}$ be a
Toeplitz-Cuntz-Krieger $E$-family such that $q_v \not= 0$ for
all $v$, and let $B := C^*(\{q_v, t_e : v \in E^0, e \in
E^1\})$. The universal property of $\Tt C^*(E)$ implies that
there is a homomorphism $\pi_{q,t} : \Tt C^*(E) \to B$
satisfying $\pi_{q,t}(\pT{v}) = q_v$ for all $v \in E^0$ and
$\pi_{q,t}(\sT{e}) = t_e$ for all $e \in E^1$. Since each $q_v$
is nonzero, $I = \ker(\pi_{q,t})$ is an ideal of $\Tt C^*(E)$
such that $\pT{v} \not\in I$ for all $v \in E^0$. Let $\kappa :
C(E) \to \TT$ and $I^\kappa$ be as in
Corollary~\ref{prop:Toeplitz big ideal}. Since $I \subset
I^\kappa$, there is a well-defined homomorphism $\psi_0 : B \to
\Tt C^*(E)/I^\kappa$ satisfying $\psi_0(q_v) = \pT{v} +
I^\kappa$ for all $v \in E^0$ and $\psi_0(t_e) = \sT{e} +
I^\kappa$ for all $e \in E^1$. Let
$\widetilde{\tau_{\overline{\kappa}}} : \Tt C^*(E)/I^1 \to \Tt
C^*(E)/I^\kappa$ be as in~\ref{lem:Ikappas same}. Then
$\psi_{q,t} := \widetilde{\tau_{\overline{\kappa}}}^{-1} \circ
\psi_0$ has the desired property. This proves
statement~\ref{it:existence}.

For statement~\ref{it:injectivity}, suppose first that
$\psi_{q,t}$ is injective. For each $\mu \in C(E)$, we have
$S_\mu = P_{r(\mu)}$ by definition of $I^1$. Let $x(\mu)$ be
the unique element of the cutting set $X$ which belongs to
$[\mu]^1$. With $\psi_{q,t}$ and $\kappa$ are as
in~(\ref{it:existence}), we have $\psi_{q,t}(t_\mu) =
\kappa(x(\mu)) \psi_{q,t}(q_{r(\mu)})$. Since $\psi_{q,t}$ is
injective, we must have $t_\mu = \kappa(x(\mu)) q_{r(\mu)}$ for
all $\mu \in C(E)$. Now suppose that there is a function
$\kappa : C(E) \to \TT$ such that $t_\mu = \kappa([\mu])
q_{r(\mu)}$ for every cycle $\mu$ with no entrance in $E$. Then
the kernel $I_{q,t}$ of the canonical homomorphism $\pi_{q,t} :
\Tt C^*(E) \to B$ contains the generators of $I^\kappa$, and
hence contains $I^\kappa$. Since the $q_v$ are all nonzero,
Corollary~\ref{prop:Toeplitz big ideal} implies that we also
have $I_{q,t} \subset I^{\lambda}$ for some $\lambda : C(E) \to
\TT$. We claim that $\kappa = \lambda$; for if not, then there
exists $C \in C(E)$ such that $k := \kappa(C)$ is distinct from
$l = \lambda(C)$. For $\mu \in C$, we then have $k\pT{v} -
\sT{\mu}, l\pT{v} - \sT{\mu} \in I^\lambda$. But then $(k - l)
\pT{v} \in I^\lambda$, which is impossible by
Lemma~\ref{lem:Ikappa vert proj free}. Hence $I^\lambda =
I_{q,t} = I^\kappa$, and $\psi_{q,t}$ is injective.
\end{proof}

\section{Properties of the co-universal
\texorpdfstring{$C^*$}{C*}-algebra}\label{sec:properties}

In this section we prove a uniqueness theorem for $\Cr{E}$ in
terms of its co-universal property. We go on to explore the
structure and properties of the co-universal algebra.
Throughout this section we have preferred proofs which
emphasise the utility of the co-universal property over other
techniques.

Let $E$ be a row-finite directed graph. We say that a
Cuntz-Krieger $E$-family $\{p_v : v \in E^0\}$, $\{s_e : e \in
E^1\}$ is a \emph{reduced Cuntz-Krieger $E$-family} if
\begin{itemize}
\item[(R)] for every cycle $\mu$ without an entrance in
    $E^0$, there is a scalar $\kappa(\mu) \in \TT$ such
    that $s_\mu = \kappa(\mu)p_{r(\mu)}$.
\end{itemize}
We say that $\{p_v : v \in E^0\}$, $\{s_e : e \in E^1\}$ is a
\emph{normalised reduced Cuntz-Krieger $E$-family} if $s_\mu =
p_{r(\mu)}$ for each cycle $\mu$ without an entrance in $E^0$.

\begin{theorem}\label{thm:Cr(E)-Properties}
Let $E$ be a row-finite directed graph.
\begin{enumerate}
\item\label{it:generation} There is a normalised reduced
    Cuntz-Krieger $E$-family $\{\pr{v} : v \in E^0\} \cup
    \{\sr{e} : E \in E^1\}$ which generates $\Cr{E}$ and
    satisfies
    Theorem~\ref{thm:Cr(E)-Existence}(\ref{it:generators})~and~(\ref{it:co-universal}).
    In particular, given any cutting set $X$ for $E$,
    $\Cr{E}$ is generated by $\{\pr{v} : v \in E^0\} \cup
    \{\sr{e} : E \in E^1 \setminus X\}$.
\item\label{it:uniqueness} Any other $C^*$-algebra
    generated by a Toeplitz-Cuntz-Krieger $E$-family
    satisfying
    Theorem~\ref{thm:Cr(E)-Existence}(\ref{it:generators})~and~(\ref{it:co-universal})
    is isomorphic to $\Cr{E}$.
\end{enumerate}
\end{theorem}
\begin{proof}
For~(\ref{it:generation}) let $\{P^\infty_v : v \in E^0\}$,
$\{S^\infty_e : e \in E^1\}$ be the Cuntz-Krieger $E$-family of
Notation~\ref{ntn:bdryCKfam}.
Theorem~\ref{thm:Cr(E)-Existence}(\ref{it:co-universal})
ensures that there is a function $\kappa : X \to \TT$ and a
homomorphism $\psi_{P^\infty, S^\infty}$ from
$C^*(\{P^\infty_v, S^\infty_e : v \in E^0, e \in E^1\})$ onto
$\Cr{E}$ such that $\psi_{P^\infty, S^\infty}(P^\infty_v) =
P_v$ for all $v \in E^0$, $\psi_{P^\infty,
S^\infty}(S^\infty_e) = S_e$ for all $e \in E^1 \setminus X$,
and $\psi_{P^\infty, S^\infty}(S^\infty_x) = \kappa(x)S_x$ for
all $x \in X$. Hence $\pr{v} := \psi_{P^\infty,
S^\infty}(P^\infty_v)$ and $\sr{e} := \psi_{P^\infty,
S^\infty}(S^\infty_e)$ generate $\Cr{E}$ and satisfy
Theorem~\ref{thm:Cr(E)-Existence}(\ref{it:generators})~and~(\ref{it:co-universal}).
To prove the last assertion of~(\ref{it:generation}), fix $x
\in X$ and calculate:
\begin{align*}
S^\infty_x
    &= S^\infty_x S^\infty_{\lambda(x)} (S^\infty_{\lambda(x)})^* \\
    &= S^\infty_{\mu(x)} (S^\infty_{\lambda(x)})^* \\
    &= (S^\infty_{\lambda(x)})^* \\
    &\in C^*(\{P^\infty_v, S^\infty_e : v \in E^0, E \in E^1 \setminus X\}.
\end{align*}

For~(\ref{it:uniqueness}), let $A$ be another $C^*$-algebra
generated by a Toeplitz-Cuntz-Krieger family $\{p^A_v : v \in
E^0\}$, $\{s^A_e : e \in E^1\}$ with each $p^A_v$ nonzero, and
suppose that $A$ has the same two properties as $\Cr{E}$.
Applying the co-universal properties, we see that there are
surjective homomorphisms $\phi : \Cr{E} \to A$ and $\psi : A
\to \Cr{E}$ such that $\phi(\pr{v}) = p^A_v$, $\phi(\sr{e}) =
s^A_e$, $\psi(p^A_v) = \pr{v}$, and $\psi(s^A_e) = \sr{e}$ for
all $v \in E^0$ and $e \in E^1 \setminus X$. In particular,
$\phi$ and $\psi$ are inverse to each other, and hence are
isomorphisms.
\end{proof}

\begin{rmk}
Of course statement~(\ref{it:generation}) of
Theorem~\ref{thm:Cr(E)-Properties} follows from the definition
of $\Cr{E}$ (embedded in the proof of
Theorem~\ref{thm:Cr(E)-Existence}). However the argument given
highlights how it follows from the co-universal property.
\end{rmk}

\begin{cor}\label{cor:injectivity criterion}
Let $E$ be a row-finite directed graph. Let $\phi : \Cr{E} \to
B$ be a homomorphism. Then $\phi$ is injective if and only if
$\phi(\pr{v}) \not= 0$ for all $v \in E^0$.
\end{cor}
\begin{proof}
Suppose that $\phi$ is injective. Then that each $\pr{v} \not=
0$ implies that each $\phi(\pr{v}) \not= 0$ also.

Now suppose that $\phi(\pr{v}) \not= 0$ for all $v \in E^0$.
Then the co-universal property of $\Cr{E}$ ensures that for any
cutting set $X$ for $E$, there is a homomorphism $\psi :
\phi(\Cr{E}) \to \Cr{E}$ satisfying $\psi(\phi(\pr{v})) =
\pr{v}$ for all $v \in E^0$ and $\psi(\phi(\sr{e})) = \sr{e}$
for all $e \in E^1 \setminus X$.
Theorem~\ref{thm:Cr(E)-Properties}(\ref{it:generation}) implies
that $\psi$ is surjective and an inverse for $\phi$.
\end{proof}

\begin{cor}\label{cor:simplicity}
Let $E$ be a row-finite directed graph. The $C^*$-algebra
$\Cr{E}$ is simple if and only if $E$ is cofinal.
\end{cor}
\begin{proof}
First suppose that $E$ is cofinal. Fix a homomorphism $\phi :
\Cr{E} \to B$. We must show that $\phi$ is either trivial or
injective. The argument of \cite[Proposition~4.2]{Raeburn2005}
shows that $\phi(\pr{v}) = 0$ for any $v \in E^0$, then
$\phi(\pr{w}) = 0$ for all $v \in E^0$, which forces $\phi =
0$. On the other hand, if $\phi(\pr{w}) \not = 0$ for all $w
\in E^0$, then Corollary~\ref{cor:injectivity criterion}
implies that $\phi$ is injective.

Now suppose that $E$ is not cofinal. Fix $v \in E^0$ and $x \in
E^{\le \infty}$ such that $v E^* x(n) = \emptyset$ for all $n
\in \NN$. Standard calculations show that
\[
I_x := \clsp\{\sr{\alpha} \srstar{\beta} : s(\alpha) = s(\beta) = x(n)\text{ for some } n \le |x|\}
\]
is an ideal of $\Cr{E}$ which is nontrivial because it contains
$\pr{x(0)}$. To see that $\pr{v} \not\in I_x$, fix $n \le |x|$
and $\alpha, \beta \in E^* x(n)$. It suffices to show that
$\pr{v} \sr{\alpha} \srstar{\beta} = 0$. Let $l := |\alpha|$.
By Theorem~\ref{thm:Cr(E)-Existence}(\ref{it:injectivity}),
$\{\pr{v} : v \in E^0\}$, $\{\sr{e}, e \in E^1\}$ is a reduced
Cuntz-Krieger $E$-family, and in particular a standard
inductive argument based on relation~(CK) shows that
\[
\pr{v} = \sum_{\lambda \in v E^{\le l}} \sr{\lambda}\srstar{\lambda}.
\]
Fix $\lambda \in v E^{\le l}$. Since $v E^* x(n) = \emptyset$
for all $n \in \NN$, we have $\alpha \not= \lambda\lambda'$ for
all $\lambda' \in E^*$. Since $|\alpha| = l \ge |\lambda|$, it
follows that $\srstar{\lambda} \sr{\alpha} = 0$. Hence $\pr{v}
\sr{\alpha} \srstar{\beta} = 0$ as claimed.
\end{proof}

\begin{cor}\label{cor:universal property}
Let $E$ be a row-finite directed graph. The $C^*$-algebra
$\Cr{E}$ is the universal $C^*$-algebra generated by a
normalised reduced Cuntz-Krieger $E$-family. That is, if $\{q_v
: v \in E^0\}$, $\{t_e : e \in E^1\}$ is another normalised
reduced Cuntz-Krieger family in a $C^*$-algebra $B$, then there
is a homomorphism $\pi_{q,t} : \Cr{E} \to B$ such that
$\pi_{q,t}(\pr{v}) = q_v$ for all $v \in E^0$ and
$\pi_{q,t}(\sr{e}) = t_e$ for all $e \in E^1$.
\end{cor}
\begin{proof}
The universal property of $\Tt C^*(E)$ implies that there is a
homomorphism $\pi^\Tt_{q,t} : \Tt C^*(E) \to B$ such that
$\pi^\Tt_{q,t}(\pT{v}) = q_v$ and $\pi^\Tt_{q,t}(\sT{e}) = t_e$
for all $v \in E^0$ and $e \in E^1$. Let $I_{q,t} :=
\ker(\pi^\Tt_{q,t})$, and let $I_{\pr{},\sr{}}$ be the kernel
of the canonical homomorphism $\pi^\Tt_{\pr{}, \sr{}} : \Tt
C^*(E) \to \Cr{E}$. Let $K := I_{q,t} \cap I_{\pr{},\sr{}}$.
Define $p^K_v := \pT{v} + K$ and $s^K_e := \sT{e} + K$ for all
$v \in E^0$ and $e \in E^1$. Since both $\{\pr{v}, \sr{e}\}$
and $\{q_v, t_e\}$ are normalised reduced Cuntz-Krieger
families, $\{p^K_v, s^K_e\}$ is also. Since no $\pT{v}$ belongs
to $I_{\pr{}, \sr{}}$, each $p^K_v$ is nonzero. Hence
Theorem~\ref{thm:Cr(E)-Existence}(\ref{it:co-universal})~and~(\ref{it:injectivity})
imply that there is an isomorphism $\psi_{p^K, s^K} : \Tt
C^*(E)/K \to \Cr{E}$ such that $\psi_{p^K, s^K}(p^K_v) =
\pr{v}$ and $\psi_{p^K, s^K}(s^K_e) = \sr{e}$ for all $v, e$.
By definition of $K$, the homomorphism $\pi^\Tt_{q,t} : \Tt
C^*(E) \to B$ descends to a homomorphism
$\widetilde{\pi^\Tt_{q,t}} : \Tt C^*(E)/K \to B$, and then
$\pi_{q,t} := \pi^\Tt_{q,t} \circ (\psi_{p^K, s^K})^{-1}$ is
the desired homomorphism.
\end{proof}

\begin{lemma}\label{lem:EX-fam -> E-fam}
Let $E$ be a row-finite directed graph. Fix a cutting set $X$
for $E$. Define a directed graph $F$ as follows:
\begin{gather*}
F^0 = \{\zeta(v) : v \in E^0\} \\
F^1 = \{\zeta(e) : e \in E^1 \setminus X\} \\
s(\zeta(e)) = \zeta(s(e))\quad\text{and}\quad r(\zeta(e)) =
\zeta(r(e)).
\end{gather*}
There is an isomorphism $\phi$ from $C^*(F)$ to $\Cr{E}$ such
that $\phi(p_{\zeta(v)}) = \pr{v}$ for all $v \in E^0$ and
$\phi(s_{\zeta(e)}) = \sr{e}$ for all $e \in F^1$.
\end{lemma}
\begin{proof}
Let $\{p_{\zeta(v)} : v \in E^0\}$, $\{s_{\zeta(e)} : e \in E^1
\setminus X\}$ denote the universal generating Cuntz-Krieger
$F$-family in $C^*(F)$. Recall that for $x \in X$, we write
$\mu(x)$ for the unique cycle with no entrance in $E$ such that
$\mu(x)_1 = x$, and we define $\lambda(x)$ to be the path such
that $\mu(x) = x\lambda(x)$. For $\nu \in E^*$ with $|\nu| \ge
2$ and $\nu_i \not\in X$ for all $i$, we write $\zeta(\nu)$ for
the path $\zeta(\nu_1)\cdots\zeta(\nu_{|\nu|}) \in F$. Define
\begin{align*}
q_v &:= p_{\zeta(v)}\text{ for all }v \in E^0,\\
t_e &:= s_{\zeta(e)}\text{ for all }e \in E^1 \setminus X,\text{ and}\\
t_x &:= s_{\zeta(\lambda(x))}^*\text{ for all }x \in X.
\end{align*}
It suffices to show that the $q_v$ and $t_e$ form a normalised
reduced Cuntz-Krieger $E$-family; the result will then follow
from
Theorem~\ref{thm:Cr(E)-Existence}(\ref{it:co-universal})~and~(\ref{it:injectivity}).

The $q_v$ are mutually orthogonal projections because the
$p_{\zeta(v)}$ are. This establishes~(\ref{it:TR1}).

For $e \in F^1$ we have $t^*_e t_e = s^*_{\zeta(e)}
s_{\zeta(e)} = p_{s(\zeta(e))} = q_{s(e)}$. For each $x \in X$,
since $\mu(x)$ has no entrance in $E$, we have
$r(x)(F)^{|\lambda(x)|} = \{\zeta(\lambda(x))\}$, so the
Cuntz-Krieger relation forces $s_{\zeta(\lambda(x))}
s_{\zeta(\lambda(x))}^* = p_{\zeta(s(x))}$. Hence
\[
t^*_x t_x
    = s_{\zeta(\lambda(x))} s_{\zeta(\lambda(x))}^*
    = p_{\zeta(s(x))} = q_{s(x)}.
\]
This establishes~(\ref{it:TR2}).

Fix $v \in F^0$ such that $vE^1 \not= \emptyset$. If $v = r(x)$
for some $x \in X$, then $r_E^{-1}(v) = \{x\}$, and we have
\[
q_v = p_{\zeta(v)}
    = s^*_{\zeta(\lambda(x))} s_{\zeta(\lambda(x))}
    = t_x t^*_x
    = \sum_{e \in r_E^{-1}(v)} t_e t^*_e.
\]
If $v \not= r(x)$ for all $x \in X$, then $vF^1 = \{\zeta(e) :
e \in vE^1\}$, and so
\[
q_v = p_{\zeta(v)}
    = \sum_{f \in vF^1} s_f s^*_f
    = \sum_{e \in vE^1} t_e t^*_e.
\]
This establishes both (\ref{it:TR3})~and~(CK).
\end{proof}

\begin{cor}\label{cor:bdry rep}
Let $E$ be a row-finite directed graph. There is an isomorphism
\[
\psi_{P^\infty, S^\infty} : \Cr{E} \to C^*(\{P^\infty_v, S^\infty_e : v \in E^0, e \in E^1\})
\]
satisfying $\psi_{P^\infty, S^\infty}(\pr{v}) = P^\infty_v$ for
all $v \in E^0$ and $\psi_{P^\infty, S^\infty}(\sr{e}) =
S^\infty_e$ for all $e \in E^1$.
\end{cor}
\begin{proof}
As observed above, $\{P^\infty_v : v \in E^0\}$, $\{S^\infty_e
: e \in E^1\}$ is a normalised reduced Cuntz-Krieger $E$-family
with each $P^\infty_v$ nonzero. The result therefore follows
from Corollaries \ref{cor:injectivity
criterion}~and~\ref{cor:universal property}.
\end{proof}

We now identify a subspace of $\ell^2(E^{\le \infty})$ which is
invariant under the Cuntz-Krieger family of
Notation~\ref{ntn:bdryCKfam}. We use the resulting
Cuntz-Krieger family to construct a faithful conditional
expectation from $\Cr{E}$ onto its diagonal subalgebra.

Let $\Erinfty$ denote the collection
\begin{align*}
\Erinfty = \{&\alpha \in E^* : s(\alpha)E^1 = \emptyset\} \\
            &\cup \{\alpha\mu^\infty : \alpha \in E^*, \mu\text{ is a cycle with no entrance }\} \\
            &\cup \{x \in E^\infty : x \not= \alpha\rho^\infty\text{ for any }\alpha,\rho \in E^*\text{ such that } s(\alpha) = r(\rho) = s(\rho)\}.
\end{align*}
So $x \in E^{\le \infty}$ belongs to $\Erinfty$ if and only if
either $x$ is aperiodic, or $x$ has the form $\alpha\mu^\infty$
for some cycle $\mu$ with no entrance in $E$. Observe that
\begin{equation}\label{eq:STE invariance}
\text{if $x \in \Erinfty$ and if $y \in E^{\le\infty}$ and $m,n \in
\NN$ satisfy $\sigma^m(x) = \sigma^n(y)$, then $y \in
\Erinfty$.}
\end{equation}

We regard $\ell^2(\Erinfty)$ as a subspace of
$\ell^2(E^{\le\infty})$. The condition~\eqref{eq:STE
invariance} implies that $\ell^2(E^{\le\infty})$ is invariant
for the Cuntz-Krieger $E$-family of
Notation~\ref{ntn:bdryCKfam}. We may therefore define a
Cuntz-Krieger $E$-family $\{\Pap{v} : v \in E^0\}$, $\{\Sap{e}
: e \in E^1\}$ in $\Bb(\ell^2(\Erinfty))$ by
\[
\Pap{v} = P^\infty_v|_{\ell^2(\Erinfty)}
\qquad\text{ and }\qquad
\Sap{e} = S^\infty_e|_{\ell^2(\Erinfty)}
\]
for all $v \in E^0$ and $e \in E^1$. Since every vertex of $E$
is the range of at least one element of $\Erinfty$, we have
$\Pap{v} \not= 0$ for all $v \in E^0$.

\begin{lemma}\label{lem:abCK rep}
Let $E$ be a row-finite directed graph. There is an isomorphism
$\psi_{\Pap{}, \Sap{}} : \Cr{E} \to C^*(\{\Pap{v}, \Sap{e} : v
\in E^0, e \in E^1\})$ satisfying $\psi_{\Pap{},
\Sap{}}(\pr{v}) = \Pap{v}$ for all $v \in E^0$ and
$\psi_{\Pap{},\Sap{}}(\sr{e}) = \Sap{e}$ for all $e \in E^1$.
\end{lemma}
\begin{proof}
The proof is identical to that of Corollary~\ref{cor:bdry rep}.
\end{proof}

For the next proposition, let $\Ered$ denote the collection of
paths $\alpha \in E^*$ such that $\alpha \not= \beta\mu$ for
any $\beta \in E^*$ and any cycle $\mu$ with no entrance in
$E$.

\begin{prop}
Let $E$ be a row-finite directed graph.
\begin{enumerate}
\item\label{it:Cr spanning} The $C^*$-algebra $\Cr{E}$
    satisfies
    \[
    \Cr{E} = \clsp\{\sr{\alpha} \srstar{\beta} : \alpha,\beta \in \Ered, s(\alpha) = s(\beta)\}.
    \]
\item\label{it:FCE} Let $D :=
    \clsp\{\sr{\alpha}\srstar{\alpha}
    : \alpha \in E^*\}$. There is a faithful conditional expectation
    $\Psi : \Cr{E} \to D$ such that
    \[
    \Psi(\sr{\alpha} \srstar{\beta}) =
    \begin{cases}
        \sr{\alpha}\srstar{\alpha} &\text{ if $\alpha = \beta$} \\
        0 &\text{ otherwise}
    \end{cases}
    \]
    for all $\alpha,\beta \in \Ered$ with $s(\alpha) =
    s(\beta)$.
\end{enumerate}
\end{prop}
\begin{proof}
By Lemma~\ref{lem:abCK rep} it suffices to prove the
corresponding statements for the $C^*$-algebra $B :=
C^*(\{\Pap{v}, \Sap{e} : v \in E^0, e \in E^1\}$.

(\ref{it:Cr spanning}) We have $B =
\clsp\{\Sap{\alpha}\Sapstar{\beta} : \alpha,\beta \in E^*\}$
because the same is true of $\Tt C^*(E)$. If $\alpha \in E^*
\setminus \Ered$, then $\alpha = \alpha'\mu^n$ for some
$\alpha' \in \Ered$, some cycle $\mu$ with no entrance in $E$
and some $n \in \NN$. Since $\{\Pap{v} : v \in E^0\}$,
$\{\Sap{e} : e \in E^1\}$ is a normalised reduced Cuntz-Krieger
$E$-family, $\Sapsp{\mu}{n} = \Pap{r(\mu)}$, so $\Sap{\alpha} =
\Sap{\alpha'}$.

(\ref{it:FCE}) Let $\{\xi_x : x \in \Erinfty\}$ denote the
standard orthonormal basis for $\ell^2(\Erinfty)$. For each $x
\in \Erinfty$, let $\theta_{x,x} \in \Bb(\ell^2(\Erinfty))$
denote the rank-one projection onto $\CC\xi_x$. Let $\Psi$
denote the faithful conditional expectation on
$\Bb(\ell^2(\Erinfty))$ determined by $\Psi(T) = \sum_{x \in
\Erinfty} \theta_{x,x} T \theta_{x,x}$, where the convergence
is in the strong operator topology. It suffices to show that
\begin{equation}\label{eq:Psi eqn}
    \Psi(\Sap{\alpha} \Sapstar{\beta}) =
    \begin{cases}
        \Sap{\alpha}\Sapsp{\alpha}{*} &\text{ if $\alpha = \beta$} \\
        0 &\text{ otherwise}
    \end{cases}
\end{equation}
for all $\alpha,\beta \in \Ered$ with $s(\alpha) = s(\beta)$.

Fix $\alpha,\beta \in \Ered$ with $s(\alpha) = s(\beta)$. If
$\alpha = \beta$, then
\[
\Sap{\alpha}\Sapsp{\alpha}{*} = \sum_{y \in
s(\alpha)\Erinfty} p_{\alpha x},
\]
and~\eqref{eq:Psi eqn} is immediate. So suppose that $\alpha
\not= \beta$. For $x \in \Erinfty$, we have
\[
\theta_{x,x} \Sap{\alpha}\Sapsp{\beta}{*} \theta_{x,x} =
\begin{cases}
    \theta_{x,x} &\text{ if $x = \alpha y = \beta y$} \\
    0 &\text{ otherwise.}
\end{cases}
\]
Hence we must show that $\alpha y \not= \beta y$ for all $y \in
s(\alpha) \Erinfty$. Fix $y \in s(\alpha) \Erinfty$. First
observe that if $|\alpha| = |\beta| = l$, then $(\alpha y)(0,
l) = \alpha \not= \beta = (\beta y)(0,l)$. Now suppose that
$|\alpha| \not= |\beta|$; we may assume without loss of
generality that $|\alpha| < |\beta|$. We suppose that $\alpha y
= \beta y$ and seek a contradiction. That $\alpha y = \beta y$
implies that $\beta = \alpha\beta'$ and $y = \beta' y$. Hence
$r(\beta') = s(\beta')$ and $y = (\beta')^\infty$. Since $y \in
\Erinfty$, it follows that $\beta' = \mu^n$ for some cycle
$\mu$ with no entrance and some $n \in \NN$, contradicting
$\beta \in \Ered$.
\end{proof}

\end{document}